\documentclass{amsart}
\usepackage{amssymb}
\usepackage{amsmath}
\usepackage{amsfonts}

\newtheorem{theorem}{Theorem}
\theoremstyle{plain}

\newtheorem{corollary}{Corollary}

\numberwithin{equation}{section}
\input{tcilatex}

\begin{document}
\title[\textbf{a }$p$\textbf{-adic analogue of...}]{$p$\textbf{-adic
interpolation function related to multiple generalized Genocchi numbers}}
\author[\textbf{S. Araci}]{\textbf{Serkan Araci}}
\address{\textbf{University of Gaziantep, Faculty of Science and Arts,
Department of Mathematics, 27310 Gaziantep, TURKEY}}
\email{\textbf{mtsrkn@hotmail.com; mtsrkn@gmail.com; saraci88@yahoo.com.tr}}
\author[\textbf{M. Acikgoz}]{\textbf{Mehmet Acikgoz}}
\address{\textbf{University of Gaziantep, Faculty of Science and Arts,
Department of Mathematics, 27310 Gaziantep, TURKEY}}
\email{\textbf{acikgoz@gantep.edu.tr}}
\author[\textbf{E. \c{S}en}]{\textbf{Erdo\u{g}an \c{S}en}}
\address{\textbf{Department of Mathematics, Faculty of Science and Letters,
Namik Kemal University, 59030 Tekirda\u{g}, TURKEY}}
\email{\textbf{erdogan.math@gmail.com}}

\begin{abstract}
In the present paper, we deal with multiple generalized Genocchi numbers and
polynomials. Also, we introduce analytic interpolating function for the
multiple generalized Genocchi numbers attached to $\chi $ at negative
integers in complex plane and we define the multiple Genocchi $p$-adic $L$%
-function. Finally, we derive the value of the partial derivative of our
multiple $p$-adic $l$-function at $s=0$.

\vspace{2mm}\noindent \textsc{2010 Mathematics Subject Classification.}
11S80, 11B68.

\vspace{2mm}

\end{abstract}

\maketitle

\section{\textbf{PRELIMINARIES}}

The works of generalized Bernoulli, Euler and Genocchi numbers and
polynomials and their combinatorial relations have received much attention 
\cite{Araci 6}, \cite{Kim 1}-\cite{Kim 25}, \cite{Liu}, \cite{Frappier}, 
\cite{Jolany 1}, \cite{Jolany 2}, \cite{Vandiver}. Generalized Bernoulli
polynomials, generalized Euler polynomials and generalized Genocchi numbers
and polynomials are the signs of very strong relationship between elementary
number theory, complex analytic number theory, Homotopy theory (stable
Homotopy groups of spheres), differential topology (differential structures
on spheres), theory of modular forms (Eisenstein series), $p$-adic analytic
numbers theory ($p$-adic $L$-functions), quantum physics(quantum Groups).

$p$-adic numbers also were invented by German Mathematician Kurt Hensel
around the end of the nineteenth century. In spite of their being already
one hundred years old, these numbers are still today enveloped in an aura of
mystery within the scientific community. The $p$-adic integral was used in
mathematical physics, for instance, the functional equation of the $q$-zeta
function, $q$-stirling numbers and $q$-Mahler theory of integration with
respect to the ring $%
\mathbb{Z}
_{p}$ together with Iwasawa's $p$-adic $L$ functions.

Also the $p$-adic interpolation functions of the Bernoulli and Euler
polynomials have been treated by Tsumura \cite{Tsumura} and Young \cite%
{Young}. T. Kim \cite{Kim 1}-\cite{Kim 24} also studied on $p$-adic
interpolation functions of these numbers and polynomials. In \cite{Carlitz},
Carlitz originally constructed $q$-Bernoulli numbers and polynomials. These
numbers and polynomials are studied by many authors (see cf. \cite{kim 2}-%
\cite{Park}, \cite{Srivastava 1}). In the last decade, a surprising number
of papers appeared proposing new generalizations of the Bernoulli, Euler and
Genocchi polynomials to real and complex variables.

In \cite{Kim 1}-\cite{Ryoo}, Kim studied some families of multiple
Bernoulli, Euler and Genocchi numbers and polynomials. By using the
fermionic $p$-adic invariant integral on $%
\mathbb{Z}
_{p}$, he constructed $p$-adic Bernoulli, $p$-adic Euler and $p$-adic
Genocchi numbers and polynomials of higher order.

In this paper, by using Kim's method in \cite{Kim 21}, we derive several
properties for the multiple generalized Genocchi numbers attached to $\chi $.

As is well-known, Genocchi numbers are defined in the complex plane by the
following exponential function%
\begin{equation}
C\left( t\right) =\frac{2t}{e^{t}+1}=\sum_{n=0}^{\infty }G_{n}\frac{t^{n}}{n!%
}\text{, }\left\vert t\right\vert <\pi \text{.}  \label{equation 1}
\end{equation}

It follows from the description that $G_{0}=0,$ $G_{1}=1,$ $G_{2}=-1,$ $%
G_{3}=0,$ $G_{4}=1,$ $G_{5}=0,\cdots ,$ and $G_{2k+1}=0$ for $k=1,2,3,\cdots
.$

The Genocchi polynomials are also given by the rule:%
\begin{equation*}
C\left( t,x\right) =e^{tG\left( x\right) }=\sum_{n=0}^{\infty }G_{n}\left(
x\right) \frac{t^{n}}{n!}=\frac{2t}{e^{t}+1}e^{xt}\text{,}
\end{equation*}%
with the usual convention of replacing $G^{n}\left( x\right) :=G_{n}\left(
x\right) $ (see \cite{kim 2}, \cite{kim 3} and \cite{kim6}).

Let $w\in 
\mathbb{N}
$. Then the multiple Genocchi polynomials of order $w$ are given by \cite%
{Ryoo}%
\begin{equation}
C^{\left( w\right) }\left( t,x\right) =\left( \frac{2t}{e^{t}+1}\right)
^{w}e^{xt}=\sum_{n=0}^{\infty }G_{n}^{\left( w\right) }\left( x\right) \frac{%
t^{n}}{n!}\text{, }\left\vert t\right\vert <\pi \text{.}  \label{equation 2}
\end{equation}

Taking $x=0$ in (\ref{equation 2}), then we have $G_{n}^{\left( w\right)
}\left( 0\right) :=G_{n}^{\left( w\right) }$ are called the multiple
Genocchi numbers of order $w$.

For $f\in 
\mathbb{N}
$ with $f\equiv 1\left( \func{mod}2\right) $, we assume that $\chi $ is a
primitive Dirichlet's charachter with conductor $f$. It is known in \cite%
{Jang 1} that the Genocchi numbers associated with $\chi $, $G_{n,\chi }$,
was introduced by the following expression%
\begin{equation}
C_{\chi }\left( t\right) =2t\sum_{\xi =1}^{f}\frac{\chi \left( \xi \right)
\left( -1\right) ^{\xi }e^{\xi t}}{e^{ft}+1}=\sum_{n=0}^{\infty }G_{n,\chi }%
\frac{t^{n}}{n!},\text{ }\left\vert t\right\vert <\frac{\pi }{f}\text{.}
\label{equation 3}
\end{equation}

In this paper, we contemplate the definition of the generating function of
the multiple generalized Genocchi numbers attached to $\chi $ in the complex
plane. From this definition, we introduce an analytic interpolating function
for the multiple generalized Genocchi numbers attached to $\chi $. Finally,
we investigate behaviour of analytic interpolating function at $s=0$.

\section{\textbf{ON AN ANALYTIC FUNCTION IN CONNECTION WITH THE MULTIPLE
GENERALIZED GENOCCHI NUMBERS}}

In this part, we introduce the multiple generalized Genocci numbers attached
to $\chi $ defined by 
\begin{align}
C_{\chi }^{\left( w\right) }\left( t\right) & =\sum_{n=0}^{\infty }G_{n,\chi
}^{\left( w\right) }\frac{t^{n}}{n!}  \label{equation 44} \\
& =\left( 2t\right) ^{w}\sum_{a_{1},\cdots ,a_{w}=1}^{f}\frac{\left(
-1\right) ^{a_{1}+\cdots +a_{w}}\chi \left( a_{1}+\cdots +a_{w}\right)
e^{t\left( a_{1}+\cdots +a_{w}\right) }}{\left( e^{ft}+1\right) ^{w}}\text{.}
\notag
\end{align}

On account of (\ref{equation 2}) and (\ref{equation 44}), we easily derive
the following%
\begin{equation}
G_{n,\chi }^{\left( w\right) }=\frac{f^{n}}{f^{w}}\sum_{a_{1},\cdots
,a_{w}=1}^{f}\left( -1\right) ^{a_{1}+\cdots +a_{w}}\chi \left( a_{1}+\cdots
+a_{w}\right) G_{n}^{\left( w\right) }\left( \frac{a_{1}+\cdots +a_{w}}{f}%
\right) \text{.}  \label{equation 5}
\end{equation}

For $s\in 
\mathbb{C}
$, we have%
\begin{gather}
\frac{1}{\Gamma \left( s\right) }\int_{0}^{\infty }t^{s-w-1}\left\{ \left(
-1\right) ^{w}C^{\left( w\right) }\left( -t,x\right) \right\} dt
\label{equation 6} \\
=2^{w}\sum_{n_{1},\cdots ,n_{w}\geq 0}\frac{\left( -1\right) ^{n_{1}+\cdots
+n_{w}}}{\left( x+n_{1}+\cdots +n_{w}\right) ^{s}},\text{ }x\neq
0,-1,-2,\cdots .  \notag
\end{gather}%
where $\Gamma \left( s\right) $ is Euler-Gamma function, which is defined by
the rule%
\begin{equation*}
\Gamma \left( s\right) =\int_{0}^{\infty }t^{s-1}e^{-t}dt\text{.}
\end{equation*}

Thanks to (\ref{equation 6}), we give the multiple Genocchi-zeta function as
follows: for $s\in 
\mathbb{C}
$ and $x\neq 0,-1,-2,\cdots ,$ 
\begin{equation}
\zeta _{G}^{\left( w\right) }\left( s,x\right) =2^{w}\sum_{n_{1},\cdots
,n_{w}\geq 0}\frac{\left( -1\right) ^{n_{1}+\cdots +n_{w}}}{\left(
x+n_{1}+\cdots +n_{w}\right) ^{s}}\text{.}  \label{equation 7}
\end{equation}

By (\ref{equation 2}) and (\ref{equation 6}), we see that%
\begin{equation*}
\zeta _{G}^{\left( w\right) }\left( -n,x\right) =\frac{G_{n+w}^{\left(
w\right) }\left( x\right) }{\binom{n+w}{w}w!}\text{ for }n\in 
\mathbb{N}
\text{.}
\end{equation*}

By utilizing from complex integral and (\ref{equation 44}), we obtain the
following equation: for $s\in 
\mathbb{C}
$,%
\begin{gather}
\frac{1}{\Gamma \left( s\right) }\int_{0}^{\infty }t^{s-w-1}\left\{ \left(
-1\right) ^{w}C_{\chi }^{\left( w\right) }\left( -t\right) \right\} dt
\label{equation 8} \\
=2^{w}\sum_{\underset{n_{1}+\cdots +n_{w}\neq 0}{n_{1},\cdots ,n_{w}=0}}%
\frac{\chi \left( a_{1}+\cdots +a_{w}\right) \left( -1\right) ^{n_{1}+\cdots
+n_{w}}}{\left( n_{1}+\cdots +n_{w}\right) ^{s}}\text{,}  \notag
\end{gather}%
where $\chi $ is the primitive Dirichlet's character with conductor 
\begin{equation*}
f\in 
\mathbb{N}
\text{ and }f\equiv 1\left( \func{mod}2\right) \text{.}
\end{equation*}

Because of (\ref{equation 8}), we give the definition Dirichlet's type of
multiple Genocchi $L$-function in complex plane as follows:%
\begin{equation}
L^{\left( w\right) }\left( s\mid \chi \right) =2^{w}\sum_{\underset{%
n_{1}+\cdots +n_{w}\neq 0}{n_{1},\cdots ,n_{w}=0}}^{\infty }\frac{\chi
\left( a_{1}+\cdots +a_{w}\right) \left( -1\right) ^{n_{1}+\cdots +n_{w}}}{%
\left( n_{1}+\cdots +n_{w}\right) ^{s}}\text{.}  \label{equation 9}
\end{equation}

Via the (\ref{equation 44}) and (\ref{equation 9}), we derive the following
theorem:

\begin{theorem}
For any $n\in 
\mathbb{N}
$, then we have%
\begin{equation}
L^{\left( w\right) }\left( -n\mid \chi \right) =\frac{G_{n+w,\chi }^{\left(
w\right) }\left( x\right) }{\binom{n+w}{w}w!}\text{.}  \label{equation 10}
\end{equation}
\end{theorem}

Let $s$ be a complex variable, and let $a$ and $b$ be integer with $0<a<F$
and $F\equiv 1\left( \func{mod}2\right) $.

Thus, we can consider the partial zeta function $S^{\left( w\right) }\left(
s;a_{1},\cdots ,a_{w}\mid F\right) $ as follows:%
\begin{gather}
S^{\left( w\right) }\left( s;a_{1},\cdots ,a_{w}\mid F\right)
\label{equation 11} \\
=2^{w}\sum_{\underset{m_{i}\equiv a_{i}\left( \func{mod}F\right) }{%
m_{1},\cdots ,m_{w}>0}}\frac{\left( -1\right) ^{m_{1}+\cdots +m_{w}}}{\left(
m_{1}+\cdots +m_{w}\right) ^{s}}  \notag \\
=\left( -1\right) ^{a_{1}+\cdots +a_{w}}F^{-s}\zeta _{G}^{\left( w\right)
}\left( s,\frac{a_{1}+\cdots +a_{w}}{F}\right) \text{.}  \notag
\end{gather}

\begin{theorem}
\bigskip The following identity holds:%
\begin{equation*}
S^{\left( w\right) }\left( s;a_{1},\cdots ,a_{w}\mid F\right) =\left(
-1\right) ^{a_{1}+\cdots +a_{w}}F^{-s}\zeta _{G}^{\left( w\right) }\left( s,%
\frac{a_{1}+\cdots +a_{w}}{F}\right) \text{.}
\end{equation*}
\end{theorem}

Then Dirichlet's type of multiple $L$-function can be expressed as the sum:
for $s\in 
\mathbb{C}
$%
\begin{equation}
L^{\left( w\right) }\left( s\mid \chi \right) =\sum_{a_{1},\cdots
,a_{w}=1}^{F}\chi \left( a_{1}+\cdots +a_{w}\right) S^{\left( w\right)
}\left( s;a_{1},\cdots ,a_{w}\mid F\right) \text{.}  \label{equation 12}
\end{equation}

Substituting $s=w-n$ into (\ref{equation 11}), we readily derive the
following: for $w,n\in 
\mathbb{N}
$%
\begin{eqnarray}
&&\binom{n}{w}w!S^{\left( w\right) }\left( w-n;a_{1},\cdots ,a_{w}\mid
F\right)   \label{equation 13} \\
&=&\left( -1\right) ^{a_{1}+\cdots +a_{w}}F^{n-w}G_{n}^{\left( w\right)
}\left( \frac{a_{1}+\cdots +a_{w}}{F}\right) \text{.}  \notag
\end{eqnarray}

By (\ref{equation 2}), it is simple to indicate the following%
\begin{equation}
G_{n}^{\left( w\right) }\left( x\right) =\sum_{k=0}^{n}\binom{n}{k}%
x^{n-k}G_{k}^{\left( w\right) }=\sum_{k=0}^{n}\binom{n}{k}%
x^{k}G_{n-k}^{\left( w\right) }\text{.}  \label{equation 14}
\end{equation}

Thanks to (\ref{equation 11}), (\ref{equation 13}) and (\ref{equation 14}),
we develop the following theorem:

\begin{theorem}
The following identity%
\begin{gather}
w!\binom{-s}{w}S^{\left( w\right) }\left( s+w;a_{1},\cdots ,a_{w}\mid
F\right)  \label{equation 15} \\
=\left( -1\right) ^{a_{1}+\cdots +a_{w}}F^{-w}\left( a_{1}+\cdots
+a_{w}\right) ^{-s}\sum_{k\geq 0}\binom{-s}{k}\left( \frac{F}{a_{1}+\cdots
+a_{w}}\right) ^{k}G_{k}^{\left( w\right) }  \notag
\end{gather}%
is true.
\end{theorem}

From (\ref{equation 12}), (\ref{equation 13}) and (\ref{equation 15}), we
have the following corollary:

\begin{corollary}
The following holds true:%
\begin{gather}
w!\binom{-s}{w}L^{\left( w\right) }\left( s+w\mid \chi \right)
\label{equation 16} \\
=\sum_{a_{1},\cdots ,a_{w}=1}^{F}\chi \left( a_{1}+\cdots +a_{w}\right)
\left( -1\right) ^{a_{1}+\cdots +a_{w}}F^{-w}\left( a_{1}+\cdots
+a_{w}\right) ^{-s}  \notag \\
\times \sum_{k=0}^{\infty }\binom{-s}{k}\left( \frac{F}{a_{1}+\cdots +a_{w}}%
\right) ^{k}G_{k}^{\left( w\right) }\text{.}  \notag
\end{gather}
\end{corollary}

The values of $L^{\left( w\right) }\left( s\mid \chi \right) $ at negative
integers are algebraic, hence may be regarded as lying in an extension of $%
\mathbb{Q}
_{p}$. We therefore look for a $p$-adic function which agrees with $%
L^{\left( w\right) }\left( s\mid \chi \right) $ at the negative integers in
the next section.

\section{\textbf{CONCLUSION}}

In this final section, we consider $p$-adic interpolation function of the
multiple generalized Genocchi $L$-function, which interpolate Dirichlet's
type of multiple Genocchi numbers at negative integers. Firstly, Washington
constructed $p$-adic $l$-function which interpolates generalized classical
Bernoulli numbers.

Here, we use some the following notations, which will be useful in remainder
of paper.

Let $\omega $ denote the $Teichm\ddot{u}ller$ character by the conductor $%
f_{\omega }=p$. For an arbitrary character $\chi $, we set $\chi _{n}=\chi
\omega ^{-n}$, $n\in 
\mathbb{Z}
$, in the sense of product of characters.

Let%
\begin{equation*}
\left\langle a\right\rangle =\omega ^{-1}\left( a\right) a=\frac{a}{\omega
\left( a\right) }\text{.}
\end{equation*}

Thus, we note that $\left\langle a\right\rangle \equiv 1\left( \func{mod}p%
\mathbb{Z}
_{p}\right) $. Let 
\begin{equation*}
A_{j}\left( x\right) =\sum_{n=0}^{\infty }a_{n,j}x^{n}\text{, }a_{n,j}\in 
\mathbb{C}
_{p}\text{, }j=0,1,2,...
\end{equation*}%
be a sequence of power series, each convergent on a fixed subset%
\begin{equation*}
T=\left\{ s\in 
\mathbb{C}
_{p}\mid \left\vert s\right\vert _{p}<p^{-\frac{2-p}{p-1}}\right\} \text{,}
\end{equation*}%
of $%
\mathbb{C}
_{p}$ such that

(1) $a_{n,j}\rightarrow a_{n,0}$ as\ $j\rightarrow \infty $ for any $n$;

(2) \ for each $s\in T$ and $\epsilon >0$, there exists an $%
n_{0}=n_{0}\left( s,\epsilon \right) $ such that 
\begin{equation*}
\left\vert \sum_{n\geq n_{0}}a_{n,j}s^{n}\right\vert _{p}<\epsilon \text{
for }\forall j\text{.}
\end{equation*}

So,%
\begin{equation*}
\lim_{j\rightarrow \infty }A_{j}\left( s\right) =A_{0}\left( s\right) \text{%
, for all }s\in T\text{.}
\end{equation*}

This was firstly introduced by Washington \cite{Washington} to indicate that
each functions $\omega ^{-s}\left( a\right) a^{s}$ and 
\begin{equation*}
\sum_{k=0}^{\infty }\binom{s}{k}\left( \frac{F}{a}\right) ^{k}B_{k}\text{,}
\end{equation*}%
where $F$ is multiple of $p$ and $f$ and $B_{k}$ is the $k$-th Bernoulli
numbers, is analytic on $T$ (for more information, see \cite{Washington}).

We assume that $\chi $ is a primitive Dirichlet's character\ with conductor $%
f\in 
\mathbb{N}
$ with $f\equiv 1\left( \func{mod}2\right) $. Then we contemplate the
multiple Genocchi $p$-adic $L$-function, $L_{p}^{\left( w\right) }\left(
s\mid \chi \right) $, which interpolates the multiple generalized Genocchi
numbers attached to $\chi $ at negative integers.

For $f\in 
\mathbb{N}
$ with $f\equiv 1\left( \func{mod}2\right) $, let us assume that $F$ is a
positive integral multiple of $p$ and $f=f_{\chi }$. We now give the
definition of mutiple Genocchi $p$-adic $L$-function as follows:%
\begin{gather}
w!\binom{-s}{w}L_{p}^{\left( w\right) }\left( s+w\mid \chi \right)
\label{equation 17} \\
=\frac{1}{F^{w}}\sum_{a_{1},\cdots ,a_{w}=1}^{F}\chi \left( a_{1}+\cdots
+a_{w}\right) \left( -1\right) ^{a_{1}+\cdots +a_{w}}\left\langle
a_{1}+\cdots +a_{w}\right\rangle ^{-s}  \notag \\
\times \sum_{k=0}^{\infty }\binom{-s}{k}\left( \frac{F}{a_{1}+\cdots +a_{w}}%
\right) ^{k}G_{k}^{\left( w\right) }\text{.}  \notag
\end{gather}

Due to (\ref{equation 17}), we want to note that $L_{p}^{\left( w\right)
}\left( s+w\mid \chi \right) $ is an analytic function on $s\in T$.

For $n\in 
\mathbb{N}
$, we have 
\begin{equation}
G_{n,\chi _{n}}^{\left( w\right) }=\frac{F^{n}}{F^{w}}\sum_{a_{1},\cdots
,a_{w}=1}^{F}\left( -1\right) ^{a_{1}+\cdots +a_{w}}\chi _{n}\left(
a_{1}+\cdots +a_{w}\right) G_{n}^{\left( w\right) }\left( \frac{a_{1}+\cdots
+a_{w}}{F}\right) \text{.}  \label{equation 18}
\end{equation}

If $\chi _{n}\left( p\right) \neq 0$, then $\left( p,f_{\chi _{n}}\right) =1$%
, and so the ratio $\frac{F}{p}$ is a multiple of $f_{\chi _{n}}$.

Let 
\begin{equation*}
\vartheta =\left\{ \frac{a_{1}+\cdots +a_{w}}{p}\mid a_{1}+\cdots
+a_{w}\equiv 0(\func{mod}p)\text{ for some }a_{i}\in 
\mathbb{Z}
\text{ with }0\leq a_{i}\leq F\right\} \text{.}
\end{equation*}

Therefore we can write the following%
\begin{eqnarray}
&&\frac{F^{n}}{F^{w}}\sum_{\underset{p\mid a_{1}+\cdots +a_{w}}{a_{1},\cdots
,a_{w}=1}}^{F}\left( -1\right) ^{a_{1}+\cdots +a_{w}}\chi _{n}\left(
a_{1}+\cdots +a_{w}\right) G_{n}^{\left( w\right) }\left( \frac{a_{1}+\cdots
+a_{w}}{F}\right)  \label{equation 19} \\
&=&p^{n-w}\frac{\left( \frac{F}{p}\right) ^{n}}{\left( \frac{F}{p}\right)
^{w}}\chi _{n}\left( p\right) \sum_{\underset{\lambda \in \vartheta }{%
a_{1},\cdots ,a_{w}=1}}^{\frac{F}{p}}\left( -1\right) ^{\lambda }\chi
_{n}\left( \lambda \right) G_{n}^{\left( w\right) }\left( \frac{\lambda }{F/p%
}\right) \text{.}  \notag
\end{eqnarray}

By (\ref{equation 19}), we define the different multiple generalized
Genocchi numbers attached to $\chi $ as follows:%
\begin{equation}
G_{n,\chi _{n}}^{\ast \left( w\right) }=\frac{\left( \frac{F}{p}\right) ^{n}%
}{\left( \frac{F}{p}\right) ^{w}}\sum_{\underset{\lambda \in \vartheta }{%
a_{1},\cdots ,a_{w}=1}}^{\frac{F}{p}}\left( -1\right) ^{\lambda }\chi
_{n}\left( \lambda \right) G_{n}^{\left( w\right) }\left( \frac{\lambda }{F/p%
}\right) \text{.}  \label{equation 20}
\end{equation}

On accounct of (\ref{equation 18}), (\ref{equation 19}) and (\ref{equation
20}), we attain the following%
\begin{equation}
G_{n,\chi _{n}}^{\left( w\right) }-p^{n-w}\chi _{n}\left( p\right) G_{n,\chi
_{n}}^{\ast \left( w\right) }=\frac{F^{n}}{F^{w}}\sum_{\underset{p\nshortmid
a_{1}+\cdots +a_{w}}{a_{1},\cdots ,a_{w}=1}}^{F}\left( -1\right)
^{a_{1}+\cdots +a_{w}}\chi _{n}\left( a_{1}+\cdots +a_{w}\right)
G_{n}^{\left( w\right) }\left( \frac{a_{1}+\cdots +a_{w}}{F}\right) \text{.}
\label{equation 21}
\end{equation}

By the definition of the multiple Genocchi polynomials of order $w$, we
write the following%
\begin{equation}
G_{n}^{\left( w\right) }\left( \frac{a_{1}+\cdots +a_{w}}{F}\right)
=F^{-n}\left( a_{1}+\cdots +a_{w}\right) ^{n}\sum_{k=0}^{n}\binom{n}{k}%
\left( \frac{F}{a_{1}+\cdots +a_{w}}\right) ^{k}G_{k}^{\left( w\right) }%
\text{.}  \label{equation 22}
\end{equation}

By (\ref{equation 21}) and (\ref{equation 22}), we have%
\begin{gather}
G_{n,\chi _{n}}^{\left( w\right) }-p^{n-w}\chi _{n}\left( p\right) G_{n,\chi
_{n}}^{\ast \left( w\right) }  \label{equation 23} \\
=\frac{1}{F^{w}}\sum_{\underset{p\nshortmid a_{1}+\cdots +a_{w}}{%
a_{1},\cdots ,a_{w}=1}}^{F}\left( -1\right) ^{a_{1}+\cdots +a_{w}}\chi
_{n}\left( a_{1}+\cdots +a_{w}\right) \left( a_{1}+\cdots +a_{w}\right) ^{n}
\notag \\
\times \sum_{k=0}^{n}\binom{n}{k}\left( \frac{F}{a_{1}+\cdots +a_{w}}\right)
^{k}G_{k}^{\left( w\right) }  \notag
\end{gather}

By (\ref{equation 17}) and (\ref{equation 23}), we readily see that%
\begin{gather}
w!\binom{n}{w}L_{p}^{\left( w\right) }\left( w-n\mid \chi \right)  \notag \\
=\frac{1}{F^{w}}\sum_{\underset{p\nshortmid a_{1}+\cdots +a_{w}}{%
a_{1},\cdots ,a_{w}=1}}^{F}\left( -1\right) ^{a_{1}+\cdots +a_{w}}\chi
_{n}\left( a_{1}+\cdots +a_{w}\right) \left( a_{1}+\cdots +a_{w}\right) ^{n}
\notag \\
\times \sum_{k=0}^{n}\binom{n}{k}\left( \frac{F}{a_{1}+\cdots +a_{w}}\right)
^{k}G_{k}^{\left( w\right) }  \label{equation 24} \\
=G_{n,\chi _{n}}^{\left( w\right) }-p^{n-w}\chi _{n}\left( p\right)
G_{n,\chi _{n}}^{\ast \left( w\right) }  \notag
\end{gather}

Consequently, we arrive at the following theorem.

\begin{theorem}
The following nice identity holds true:%
\begin{gather*}
w!\binom{-s}{w}L_{p}^{\left( w\right) }\left( s+w\mid \chi \right) \\
=\frac{1}{F^{w}}\sum_{a_{1},\cdots ,a_{w}=1}^{F}\chi \left( a_{1}+\cdots
+a_{w}\right) \left( -1\right) ^{a_{1}+\cdots +a_{w}}\left\langle
a_{1}+\cdots +a_{w}\right\rangle ^{-s} \\
\times \sum_{k=0}^{\infty }\binom{-s}{k}\left( \frac{F}{a_{1}+\cdots +a_{w}}%
\right) ^{k}G_{k}^{\left( w\right) }\text{.}
\end{gather*}
\end{theorem}

Thus $L_{p}^{\left( w\right) }\left( s+w\mid \chi \right) $ is an analytic
function on $T$. Additionally, for each $n\in 
\mathbb{N}
$, we procure the following:%
\begin{equation*}
w!\binom{n}{w}L_{p}^{\left( w\right) }\left( w-n\mid \chi \right) =G_{n,\chi
_{n}}^{\left( w\right) }-p^{n-w}\chi _{n}\left( p\right) G_{n,\chi
_{n}}^{\ast \left( w\right) }\text{.}
\end{equation*}

Using Taylor expansion at $s=0$, we have%
\begin{equation}
\binom{-s}{k}=\frac{\left( -1\right) ^{k}}{k}s+\cdots \text{ if }k\geq 1%
\text{.}  \label{equation 25}
\end{equation}

Differentiating on both sides in (\ref{equation 17}), with respect to $s$ at 
$s=0$, we obtain the following corollary.

\begin{theorem}
Let $F$ be a positive integral multiple of $p$ and $f$. Then we have%
\begin{gather*}
\frac{\partial }{\partial s}\left( \binom{-s}{w}L_{p}^{\left( w\right)
}\left( s+w\mid \chi \right) \right) \mid _{s=0}=\frac{\left( -1\right) ^{w}%
}{w}L_{p}^{\left( w\right) }\left( w\mid \chi \right) \\
=\frac{1}{w!F^{w}}\sum_{\underset{\left( a_{1}+\cdots +a_{w},p\right) =1}{%
a_{1},\cdots ,a_{w}=1}}^{F}\chi \left( a_{1}+\cdots +a_{w}\right) \left(
-1\right) ^{a_{1}+\cdots +a_{w}}\left( 1-\log _{p}\left( a_{1}+\cdots
+a_{w}\right) \right) \\
+\frac{1}{w!F^{w}}\sum_{\underset{\left( a_{1}+\cdots +a_{w},p\right) =1}{%
a_{1},\cdots ,a_{w}=1}}^{F}\chi \left( a_{1}+\cdots +a_{w}\right) \left(
-1\right) ^{a_{1}+\cdots +a_{w}}\sum_{m=1}^{\infty }\frac{\left( -1\right)
^{m}}{m}\left( \frac{F}{a_{1}+\cdots +a_{w}}\right) ^{k}G_{k}^{\left(
w\right) }
\end{gather*}%
where $\log _{p}x$ is the $p$-adic logarithm.
\end{theorem}

\end{document}